\documentclass[12pt]{article}
\usepackage{graphics}
\usepackage{fullpage}
\usepackage{mathptmx}

\begin{document}

\title{Algorithmic Information Theory and\\ Foundations of Probability}

\author{Alexander Shen%
\thanks{%
LIF Marseille, CNRS \& Univ.~Aix--Marseille. On leave from IITP, RAS,
Moscow. Supported in part by NAFIT
ANR-08-EMER-008-01 grant.   E-mail: \texttt{alexander.shen@lif.univ-mrs.fr}}}

\date{}
\maketitle

\begin{abstract}
The question how and why mathematical probability theory can be
applied to the ``real world'' has been debated for centuries. We
try to survey the role of algorithmic information theory
(Kolmogorov complexity) in this debate.
\end{abstract}

\section{Probability theory paradox}

One often describes the natural sciences framework as follows: a
hypothesis is used to predict something, and the prediction is
then checked against the observed actual behavior of the system;
if there is a contradiction, the hypothesis needs to be changed.

Can we include probability theory in this framework? A statistical
hypothesis (say, the assumption of a fair coin) should be then
checked against the experimental data (results of coin tossing)
and rejected if some discrepancy is found. However, there is an
obvious problem: The fair coin assumption says that in a series
of, say, $1000$ coin tossings all the $2^{1000}$ possible
outcomes (all $2^{1000}$ bit strings of length $1000$) have the
same probability $2^{-1000}$. How can we say that some of them
contradict the assumption while other do not?

The same paradox can be explained in a different way. Consider a
casino that wants to outsource the task of card shuffling to a
special factory that produced shrink-wrapped well shuffled decks
of cards. This factory would need some quality control
department. It looks at the deck before shipping it to the
customer, blocks some ``badly shuffled'' decks and approves some
others as ``well shuffled''. But how is it possible if all $n!$
orderings of $n$ cards have the same probability?

\section{Current best practice}

Whatever the philosophers say, statisticians have to perform
their duties. Let us try to provide a description of their
current ``best practice''
(see~\cite{uspensky-semenov-shen,shen}).

\textbf{A}. \emph{How a statistical hypothesis is applied}.
First of all, we have to admit that probability theory makes no
predictions but only gives recommendations: \emph{if the
probability} (computed on the basis of the statistical
hypothesis) \emph{of an event $A$ is much smaller than the
probability of an event~$B$, then the possibility of the
event~$B$ must be taken into consideration to a greater extent
than the possibility of the event~$A$} (assuming the
consequences are equally grave). For example, if the probability
of $A$ is smaller than the probability of being killed on the
street by a meteorite, we usually ignore $A$ completely (since
we have to ignore event $B$ anyway in our everyday life).

Borel~\cite{borel-1914} describes this principle as follows:
``\ldots Fewer than a million people live in Paris. Newspapers
daily inform us about the strange events or accidents that
happen to some of them. Our life would be impossible if we were
afraid of all adventures we read about. So one can say that from
a practical viewpoint we can ignore events with probability less
that one millionth\ldots\ Often trying to avoid something bad we
are confronted with even worse\ldots\ To avoid this we must know
well the probabilities of different events'' (Russian ed.,
pp.~159--160).

\textbf{B}. \emph{How a statistical hypothesis is tested}. Here
we cannot say na\"\i vely that if we observe some event that has
negligible probability according to our hypothesis, we reject
this hypothesis. Indeed, this would mean that any $1000$-bit
sequence of the outcomes would make the fair coin assumption
rejected (since this specific seqeunce has negligible
probability $2^{-1000}$).

Here algorithmic information theory comes into play: We reject
the hypothesis if we observe a \emph{simple} event that has
negligible probability according to this hypothesis. For
example, if coin tossing produces thousand tails, this event is
simple and has negligible probability, so we don't believe the
coin is fair. Both conditions (``simple'' and ``negligible
probability'') are important: the event ``the first bit is a
tail'' is simple but has probability $1/2$, so it does not
discredit the coin. On the other hand, every sequence of
outcomes has negligible probability $2^{-1000}$, but if it is
not simple, its appearance does not discredits the fair
coin assumption.

Often both parts of this scheme are united into a statement
``events with small probabilities do not happen''. For example,
Borel writes: ``One must not be afraid to use the word
``certainty'' to designate a probability that is sufficiently
close to $1$'' (\cite{borel-1950}, Russian translation, p.~7).
Sometimes this statement is called ``Cournot principle''. But we
prefer to distinguish between these two stages, because for the
hypothesis testing the existence of a simple description of an
event with negligible probability is important, and for
application of the hypothesis it seems unimportant. (We can
expect, however, that events interesting to us have simple
descriptions because of their interest.)

\section{Simple events and events specified in advance}

Unfortunately, this scheme remains not very precise: the
Kolmogorov complexity of an object $x$ (defined as the minimal
length of the program that produces $x$) depends on the choice
of programming language; we need also to fix some way to
describe the events in question. Both choices lead only to an
$O(1)$ change asymptotically; however, strictly speaking, due to
this uncertainty we cannot say that one event has smaller
complexity than the other one. (The word ``negligible'' is also
not very precise.) On the other hand, the scheme described,
while very vague, seems to be the best approximation to the
current practice.

One of the possible ways to eliminate complexity in this picture
is to say that a hypothesis is discredited if we observe a very
unprobable event \emph{that was specified in advance} (before
the experiment). Here we come to the following question. Imagine
that you make some experiment and get a sequence of thousand
bits that looks random at first. Then somebody comes and says
``Look, if we consider every third bit in this sequence, the
zeros and ones alternate''. Will you still believe in the fair
coin hypothesis? Probably not, even if you haven't thought about
this event before looking at the sequence: the event is so
simple that one \emph{could} think about it. In fact, one may
consider the union of all simple events that have small
probability, and it still has small probability (if the bound
for the complexity of a simple event is small compared to the
number of coin tossing involved, which is a reasonable condition
anyway). And this union can be considered as specified before
the experiment (e.g., in this paper).

On the other hand, if the sequence repeats some other sequence
observed earlier, we probably won't believe it is obtained by
coin tossing even if this earlier sequence had high complexity.
One may explain this opinion saying the the entire sequence of
observations is simple since it contains repetitions; however,
the first observation may be not covered by any probabilistic
assumption. This could be taked into account by considering the
\emph{conditional} complexity of the event (with respect to all
information available before the experiment).

The conclusion: we may remove one problematic requirement (being
``simple'' in a not well specified sense) and replace it by
another problematic one (being specified before the
observation).

\section{Frequency approach}

The most natural and common explanation of the notion of
probability says that probability is the limit value of frequencies
observed when the number of repetitions tends to infinity. (This
approach was advocated as the only possible basis for probability
theory by Richard von Mises.)

However, we cannot observe infinite sequences, so the actual
application of this definition should somehow deal with finite
number of repetitions. And for finite number of repetitions our
claim is not so strong: we do not guarantee that frequency of
tails for a fair coin is exactly $1/2$; we say only that it is
highly improbable that it deviates significantly from $1/2$.
Since the words ``highly improbably'' need to be interpreted,
this leads to some kind of logical circle that makes the
frequency approach much less convincing; to get out of this
logical circle we need some version of Cournot principle.

Technically, the frequency approach can be related to the
principles explained above. Indeed, the event ``the number of
tails in a $1\,000\,000$ coin tossings deviates from $500\,000$
more than by $100\,000$'' has a simple description and very
small probability, so we reject the fair coin assumption if such
an event happens (and ignore the dangers related to this event
if we accept the fair coin assumption). In this way the belief
that frequency should be close to probability (if the
statistical hypothesis is chosen correctly) can be treated as
the consequence of the principles explained above.

\section{Dynamical and statistical laws}

We have described how the probability theory is usually applied.
But the fundamental question remains: well, probability theory
describes (to some extent) the behavior of a symmetric coin or
dice and turns out to be practically useful in many cases. But
is it a new law of nature or some consequence of the known
dynamical laws of classical mechanics? Can we somehow ``prove''
that a symmetric dice indeed has the same probabilities for all
faces (if the starting point is high enough and initial linear
and rotation speeds are high enough)?

Since it is not clear what kind of ``proof'' we would like to
have, let us put the question in a more practical way. Assume
that we have a dice that is not symmetric and we know exactly
the position of its center of gravity. Can we use the laws of
mechanics to find the probabilities of different outcomes?

It seems that this is possible, at least in principle. The laws
of mecahnics determine the behavior of a dice (and therefore the
outcome) if we know the initial point in the phase space
(initial position and velocity) precisely. The phase space,
therefore, is splitted into six parts that correspond to six
outcomes. In this sense there is no uncertainty or probabilities
up to now. But these six parts are well mixed since very small
modifications affect the result, so if we consider a small (but
not very small) part of the phase space around the initial
conditions and any probability distribution on this part whose
density does not change drastically, the measures of the six
parts will follow the same proportion.

The last sentence can be transformed into a rigorous
mathematical statement if we introduce specific assumptions
about the size of the starting region in the phase space and
variations of the density of the probability distribution on it.
It then can be proved. Probably it is a rather difficult
mathematical problem not solved yet, but at least theoretically
the laws of mechanics allow us to compute the probabilities of
different outcomes for a non-symmetic dice.

\section{Are ``real'' sequences complex?}

The argument in the preceding section would not convince a
philosophically minded person. Well, we can (in principle)
compute some numbers that can be interpreted as probabilities of
the outcomes for a dice, and we do not need to fix the
distribution on the initial conditions, it is enough to assume
that this distribution is smooth enough. But still we speak
about probability distributions that are somehow externally
imposed in addition to dynamical laws.

Essentially the same question can be reformulated as follows.
Make $10^6$ coin tosses and try to compress the resulting
sequence of zeros and ones by a standard compression program,
say, gzip. (Technically, you need first to convert bit sequence
into a byte sequence.) Repeat this experiment (coin tossing plus
gzipping) as many times as you want, and this will never give
you more that $1\%$ compression. (Such a compression is possible
for less than $2^{-10000}$-fraction of all sequences.) This
statement deserves to be called a law of nature: it can be
checked experimentally in the same way as other laws are. So the
question is: does this law of nature follows from dynamical laws
we know?

To see where the problem is, it is convenient to simplify the
situation. Imagine for a while that we have discrete time, phase
space is $[0,1)$ and the dynamical law is
        $$
x \mapsto T(x)=\textbf{ if } 2x< 1 \textbf{ then } 2x \textbf{ else } 2x-1.
        $$
So we get a sequence of states $x_0, x_1=T(x_0),
x_2=T(x_1),\ldots$; at each step we observe where the current
state is~--- writing $0$ if $x_n$ is in $[0,1/2)$ and $1$ if $x_n$
is in $[1/2,1)$.

This tranformation $T$ has the mixing property we spoke about;
however, looking at it more closely, we see that a sequence of
bits obtained is just the binary representation of the initial
condition. So our process just reveals the initial condition bit
by bit, and any statement about the resulting bit sequence
(e.g., its incompressibility) is just a statement about the
initial condition.

So what? Do we need to add to the dynamical laws just one more
methaphysical law saying that world was created at the random
(=incompressible) state? Indeed, algorithmic transformations
(including dynamical laws) cannot increase significantly the
Kolmogorov complexity of the state, so if objects of high
complexity exist in the (otherwise deterministic, as we assume for now)
real world now, they should be there at
the very beginning. (Note that it is difficult to explain the
randomness observed saying that we just observe the world at
random time or in a random place: the number of bits needed to
encode the time and place in the world is not enough to explain
an incompressible string of length, say $10^6$, if we use
currently popular estimates for the size and age of the world:
the logarithms of the ratios of the maximal and minimal lengths
(or time intervals) that exist in nature are negligible compared
to $10^6$ and therefore the position in space-time cannot
determine a string of this complexity.

Should we conclude then that instead of playing the dice (as
Einstein could put it), God provided concentrated randomness
while creating the world?

\section{Randomness as ignorance: Blum -- Micali -- Yao}

This discussion becomes too philosophical to continue it
seriously. However, there is an important mathematical result
that could influence the opinion of the philosophers discussing
the notions of probability and randomness. (Unfortunately,
knowledge does not penetrate too fast, and I haven't yet seen
this argument in traditional debates about the meaning of
probability.)

This result is the existence of pseudorandom number generators
(as defined by Yao, Blum and Micali; they are standard tools in
computational cryptography, see, e.g., Goldreich
textbook~\cite{goldreich}). The existence is proved modulo some
complexity assumtions (the existence of one-way functions) that
are widely believed though not proven yet.

Let us explain what a pseudorandom number generator (in Yao --
Blum -- Micali) sense is. Here we use rather vague terms and
oversimplify the matter, but there is a rigorious
mathematics behind. So imagine a simple and fast algorithmic
procedure that gets a ``seed'', which is a binary string of
moderate size, say, $1\,000$ bits, and produces a very long
sequence of bits out of it, say, of length $10^{10}$. By
necessity the output string has small complexity compared to its
length (complexity is bounded by the seed size plus the length
of the processing program, which we assume to be rather short).
However, it may happen that the output sequences will be
``indistinguishable'' from truly random sequences of length
$10^{10}$, and in this case the transformation procedure is
called pseudorandom number generator.

It sounds as a contradiction: as we have said, output sequences
have small Kolmogorov complexity, and this property
distinguishes them from most of the sequences of length
$10^{10}$. So how they can be indistinguishable? The explanation
is that the difference becomes obvious only when we know the
seed used for producing the sequence, but there is no way to
find out this seed looking at the sequence itself. The formal
statement is quite technical, but its idea is simple. Consider
any simple test that looks at $10^{10}$-bit string and says
`yes' or `no' (by whatever reason; any simple and fast program
could be a test). Then consider two ratios: (1) the fraction of
bit strings of length $10^{10}$ that pass the test (among all bit
strings of this length); (2)~ the fraction of seeds that lead to
a $10^{10}$-bit string that passes the test (among all seeds). The
pseudorandom number generator property guarantees that these two
numbers are very close.

This implies that if some test rejects most of the pseudorandom
strings (produced by the generator), then it would also reject
most of the strings of the same length, so there is no way to
find out whether somebody gives us random or pseudorandom
strings.

In a more vague language, this example shows us that randomness
may be in the eye of the beholder, i.e., the randomness of an
observed sequence could be the consequence of our limited
computational abilities which prevent us from discovering
non-randomness. (However, if somebody shows us the seed, our
eyes are immediately opened and we see that the sequence has
very small complexity.)

In particular, trying gzip-compression on pseudorandom
sequences, we rarely would find them compressible (since
gzip-compressibility is a simple test that fails for most
sequences of length $10^{10}$, it should also fail for most
pseudorandom sequences).

So we should not exclude the possibility that the world is
governed by simple dynamical laws and its initial state can be
also described by several thousands of bits. In this case
``true'' randomness does not exist in the world, and every
sequence of $10^6$ coin tossing that happened or will happen in the
foreseeable future produces a string that has Kolmogorov
complexity much smaller than its
length. However, a computationally limited observer (like
ourselves) would never discover this fact.

\section{Digression: thermodynamics}

The connection between statistical and dynamical laws was
discussed a lot in the context of thermodynamics while
discussing the second law. However, one should be very careful
with exact definition and statements. For example, it is often
said that the Second Law of thermodynamics cannot be derived
from dynamical laws because they are time-reversible while the
second law is not. On the other hand, it is often said that the
second law has many equivalent formulations, and one of them
claims that the perpetual motion machine of the second kind is
impossible, i.e., no device can operate on a cycle to receive
heat from a single reservoir and produce a net amount of work.

However, as Nikita Markaryan explained (personal communication),
in this formulation the second law of thermodynamics \emph{is} a
consequence of dynamic laws. Here is a sketch of this argument.
Imagine a perpetual motion machine of a second kind exists.
Assume this machine is attached to a long cylinder that contains
warm gaz. Fluctuations of gaz pressure provide a heat exchange
between gaz and machine. On the other side machine has
rotating spindle and a rope to lift some weight (due to rotation).
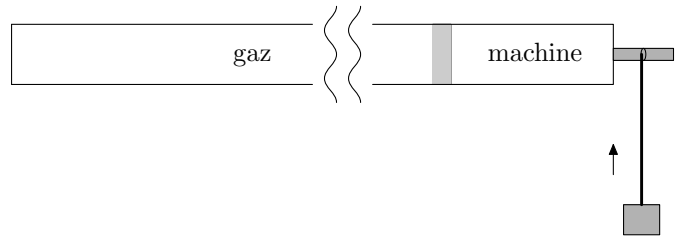
\begin{center}
        \includegraphics{machine-1.mps}
\end{center}
When the machine works, the gaz temperature (energy) goes down
and the weight goes up. This is not enough to call the machine a
perpetual motion machine of the second kind (indeed, it can
contain some amount of cold substance to cool the gaz and some
spring to lift the weight). So we assume that the rotation angle
(and height change) can be made arbitrarily large by increasing
the amount of the gaz and the length of the cylinder. We also
need to specify the initial conditions of the gaz; here the
natural requirement is that the machine works (as described) for
most initial conditions (according to the natural probability
distribibution in the gaz phase space).

Why is such a machine impossible? The phase space of the entire
system can be considered as a product of two components: the
phase space of the machine itself and the phase space of the
gaz. The components interact, and the total energy is constant.
Since the machine itself has some fixed number of components,
the dimension of its component (or the number of degrees of
freedom in the machine) is negligible compared to the dimension
of the gaz component (resp. the number of degrees of freedom in
the gaz). The phase space of the gaz is splitted into layers
corresponding to different level of energy; the higher the
energy is, the more volume is used, and this dependence
overweights the similar dependence for the machine since the gaz
has much more degrees of freedom. Since the transformation of
the phase space of the entire system is measure-preserving, it
is impossible that a trajectory started from a large set with
high probability ends in a small set: the probability of this
event does not exceed the ratio of a measures of destination and
source sets in the phase space.

This argument is quite informal and ignores many important
points. For example, the measure on the phase space of the
entire system is not exactly a product of measures on the gaz and
machine coordinates; the source set of the trajectory can have
small measure if the initial state of the machine is fixed with
very high precision, etc. (The latter case does not contradicts
the laws of thermodynamics: if the machine use a fixed amount of
cooling substance of very low temperature, the amount of work
produced can be very large.) But at least these informal
arguments make plausible that dynamic laws make imposiible the
perpetual motion machine of the second kind (if the latter is
defined properly).

\section{Digression: quantum mechanics}

Another physics topic often discussed is quantum mechanics as a
source of randomness. There were many philosophical debate
around quantum mechanics; however, it seems that the relation
between quantum mechanical models and observations resembles the
situation with probability theory and statistical mechanics; in
quantum mechanics the model assigns \emph{amplitudes} (instead of
probabilities) to different outcomes (or events). The amplitudes
are complex numbers and ``quantum Cournot principle'' says that
if the (absolute value) of the amplitude of event $A$ is smaller
than for event $B$, then the possibility of the event~$B$ must
be taken into consideration to a greater extent than the
possibility of the event~$A$ (assuming the consequences are
equally grave). Again this implies that we can (practically)
ignore events with very small amplitudes.

The interpretation of the square of amplitude as probability can
be then derived is the same way as in the case of the frequency
approach. If a system is made of $N$ independent identical
systems with two outcomes $0$ and $1$ and the outcome $1$ has
amplutude $z$ in each system, then for the entire system the
amplitude of the event ``the number of $1$'s among the outcomes
deviates significantly from $N|z|^2$'' is very small (it is just
the classical law of large numbers in disguise).

One can then try to analyze measurement devices from the quantum
mechanical viewpoint and to ``prove'' (using the same quantum
Cournot principle) that the frequency of some outcome of
measurement is close to the square of the length of the
projection of the initial state to corresponding subspace outside
some event of small amplitude, etc.

\section{Acknowledgements}
The material covered is not original: the topic was discussed
for a long time in books and papers to numerous to mention, both
the classics of the field (such as books and papers by Laplace,
Borel, and Kolmogorov) and more expository writings, such as
books written by Polya and Renyi. However, some remarks I
haven't seen before and bear responsibility for all the errors
and misunderstandings.

More information about the history of algorithmic information
theory (Kolmogorov complexity, algorithmic randomness) can be
found in~\cite{bienvenu-shen}.

I am indepted to many colleagues, including
Vladimir Fock,
Peter Gacs,
Alexey Kitaev,
Denis Kosygin,
Leonid Levin,
Leonid Levitov,
Yury Makhlin,
Nikita Mar\-ka\-ryan (who explained to me the proof of the
impossibility of the perpetuum mobile of the second kind),
Vladimir A. Uspensky,
Vladimir Vovk (who explained to me how the
complexity assumption can be replaced by the requirement that
the event should be specified before the experiment),
Michael Vyalyi,
Alexander Zvonkin,
and all the participants of Kolmogorov
seminar for many useful discussions.

I thank the organizers of LIX Colloquium ``Reachability Problems'09'' for
the opportunity to present this material at the colloquium.


\begin{thebibliography}{9}

\bibitem{bienvenu-shen}
Bienvenu, L.,  Shen, A.: Algorithmic information
theory and martingales,
preprint. Available as \texttt{arXiv:0906.2614}.

\bibitem{borel-1914}
Borel, E.: Le hazard, Alcan, Paris (1914). [Russian
translation: Gosizdat, Moscow -- Petrograd (1923)]

\bibitem{borel-1950}
Borel, E.: Probabilit\'e et
certitude, Presse Univ. de France, Paris (1950).
[English translation: Probability and certainty, Walker
publishers (1963). Russian translation: Fizmatgiz, Moscow (1961).]

\bibitem{goldreich}
Goldreich, O.: Foundations of Cryptography. Volume 1.
Basic Tools. Cambridge University Press (2007)

\bibitem{polya}
Polya, G.: Mathematics and Plausible Reasoning, vols.
1,2. Princeton University Press  (1990) [reprint]

\bibitem{renyi}
Renyi, A.: Letters on probability, Wayne State University
Press (1972)

\bibitem{uspensky-semenov-shen}
Uspensky, V.,  Semenov, A., Shen, A.: Can an
individual sequence of zeros and ones be random? Russian
Mathematical Surveys, 45(1), 121--189 (1990)

\bibitem{shen}
Shen, A.: On the logical basis of application of
probability theory, In: Proc. Workshop on semiotic aspects of
the formalization of intellectual activity, Telavi (Georgia,
then USSR), 1983, pp.~144--146. VINITI, Moscow (1983)

\end{thebibliography}
\end{document}